\documentclass[smallextended,referee,envcountsect]{svjour3}
\smartqed

\usepackage{amsmath}
\usepackage{amsfonts}
\usepackage{amssymb}
\usepackage[svgnames,dvipsnames]{xcolor}
\usepackage{bm}
\usepackage{graphicx}
\usepackage[shortlabels]{enumitem}
\setlist[enumerate]{leftmargin=.5in}
\setlist[itemize]{leftmargin=.3in}
\usepackage{booktabs}
\usepackage{mathtools}
\usepackage{algorithm}
\usepackage{algorithmicx}
\usepackage{algpseudocode}
\usepackage{float}

\usepackage[numbers,sort&compress]{natbib} 

\usepackage[nameinlink]{cleveref}
\crefname{equation}{}{} 
\crefname{thm}{Theorem}{Theorems}

\newtheorem{thm}{Theorem}[section]

\usepackage{aliascnt}
\newcommand{\newaliasthm}[4]{%
  \newaliascnt{#1}{#3}%
  \newtheorem{#1}[#1]{#2}%
  \aliascntresetthe{#1}%
  \crefname{#1}{#2}{#4}%
  \Crefname{#1}{#2}{#4}%
}

\newaliasthm{cor}{Corollary}{thm}{Corollaries}
\newaliasthm{lem}{Lemma}{thm}{Lemmas}
\newaliasthm{prop}{Proposition}{thm}{Propositions}
\newaliasthm{defn}{Definition}{thm}{Definitions}
\newaliasthm{rem}{Remark}{thm}{Remarks}

\numberwithin{equation}{section}

\DeclareMathOperator{\prox}{prox}

\DeclareMathOperator{\sign}{sign}

\newcommand{\bR}{\mathbb{R}}

\newcommand{\bS}{\mathbb{S}}

\newcommand{\vect}[1]{\bm{#1}}

\newcommand{\ve}{\vect{e}}

\newcommand{\vu}{\vect{u}}
\newcommand{\vv}{\vect{v}}

\newcommand{\vx}{\vect{x}}
\newcommand{\vy}{\vect{y}}
\newcommand{\vz}{\vect{z}}

\newcommand{\vzero}{\vect{0}}
\newcommand{\vone}{\vect{1}}

\newcommand{\cC}{\mathcal{C}}

\newcommand{\cU}{\mathcal{U}}
\newcommand{\cV}{\mathcal{V}}
\newcommand{\cX}{\mathcal{X}}

\newcommand{\mtrx}[1]{\mathsf{#1}}

\newcommand{\mI}{\mtrx{\mathrm{I}}}

\newcommand{\mP}{\mtrx{P}}
\newcommand{\mR}{\mtrx{R}}

\newcommand{\supdown}{{\scalebox{0.5}{\(\bm{\downarrow}\)}}}
\newcommand{\tvy}{\mathfrak{\vy}}
\newcommand{\ty}{\mathfrak{y}}

\newcommand{\bcC}{\widetilde{\cC}}

\newif\ifshowrevisions
\showrevisionsfalse

\journalname{JOTA}

\begin{document}


\title{Proximity Operator of the \(\ell_1\) over \(\ell_2\) Function}

\author{Lixin Shen \and Guohui Song}

\institute{Lixin Shen \at
  Department of Mathematics, Syracuse University, Syracuse, NY 13244, USA\\
  lshen03@syr.edu
  \and
  Guohui Song,  Corresponding author  \at
  Department of Mathematics and Statistics, Old Dominion University, Norfolk, VA 23529, USA\\
  gsong@odu.edu
}

\date{Received: date / Accepted: date}

\maketitle

\ifshowrevisions
  \begin{center}
    \textcolor{blue}{\small Revised material is shown in blue.}
  \end{center}
\fi

\begin{abstract}
  We study the complete, possibly set-valued proximity operator of the nonconvex, scale-invariant ratio \(h(\vx)=\|\vx\|_{1}/\|\vx\|_{2}\).
  A polar decomposition of the nonzero branch, followed by sign and permutation reductions, transforms the proximal problem into a smooth rank-one quadratic problem over the nonnegative orthant of the unit sphere.
  After the magnitudes of the data are sorted, every canonical proximal direction belongs to a finite candidate set indexed by its support size \(k\in\{1,\ldots,n\}\).
  For each \(k\geq 2\), a candidate exists if and only if an explicit quartic equation has a simple root satisfying both a positive-orthant condition and a second-order derivative condition.
  This characterization eliminates the need to guess the unknown sparsity.
  We further derive a necessary-and-sufficient \(O(1)\) prefix test that enforces positivity and prunes support sizes for which no candidate can exist.
  Using prefix sums, the optimized method performs an \(O(n)\) scan after an \(O(n\log n)\) sorting step; the construction of multiple tied outputs is necessarily output-sensitive.
  Numerical experiments illustrate the finite-candidate selection, identify inputs on which a sparsity-guessing baseline is suboptimal, quantify the benefit of pruning, and demonstrate the numerical accuracy of the secular-equation solver.
  The observed runtimes exhibit near-linear growth over the tested range, consistent with the predicted \(O(n)\) post-sorting scan and \(O(n\log n)\) total complexity.
\end{abstract}


\section{Introduction}\label{sec:intro}
The scale-invariant ratio
\begin{equation}\label{def:h}
  h(\vx)=
  \begin{cases}
    \dfrac{\|\vx\|_1}{\|\vx\|_2}, & \vx\neq\vzero, \\[1mm]
    a,                            & \vx=\vzero,
  \end{cases}
  \qquad a\in[0,1],
\end{equation}
is proper and lower semicontinuous.
Indeed, \(h\) is continuous away from the origin, and \(h(\vx)\geq1\) for \(\vx\neq\vzero\), while \(h(\vzero)=a\leq1\).
The ratio has gained traction as a surrogate for the \(\ell_0\) penalty in sparse modeling because it promotes few nonzero entries \cite{Yin-Esser-Xin:CIS14,Rahimi-Wang-Dong-Lou:SIAMSC2019} while remaining invariant under nonzero rescaling of \(\vx\).
The same ratio structure that makes \(h\) attractive also makes it nonconvex and nonsmooth, complicating both analysis and computation.

A central tool for optimization with nonsmooth terms is the proximity operator; see \cite{AttouchBolteSvaiter2013KL,Beck-Teboulle:SIAMIS:2009, BolteSabachTeboulle2014PALM,Chambolle-Pock:JMIV11,Combettes-Wajs:MMS:05,Li-Shen-Zhang-Zhou:ACHA-2022, Li-Shen-Xu-Zhang:AiCM:15,Micchelli-Shen-Xu:IP-11,Parikh-Boyd:NF-Opt:14} and the references therein.
For a proper lower semicontinuous function
\(f\), a point \(\vy\in\bR^n\), and \(\mu>0\), we use the set-valued
convention of the proximity operator
\begin{equation}\label{eq:prox}
  \prox_{\mu f}(\vy)
  :=\operatorname*{argmin}_{\vx\in\bR^n}
  Q_{f,\vy}(\vx), \qquad Q_{f,\vy}(\vx):=\frac12\|\vx-\vy\|_2^2+\mu f(\vx).
\end{equation}
A closed form or an exact low-cost procedure for this set simplifies a broad class of proximal algorithms, including primal--dual methods \cite{Chambolle-Pock:JMIV11} and FISTA \cite{Beck-Teboulle:SIAMIS:2009}.

Two lines of work have addressed \(\prox_{\mu h}\) directly.
The method in~\cite{Tao2022} derives a closed-form expression under the assumption that the \emph{true} sparsity, namely the number of nonzero entries of a proximal point, is known.
In practice, a bisection heuristic is used to estimate this sparsity, and an incorrect estimate can miss the globally optimal support.
A second approach in~\cite{Jia2024} exploits the scale and signed-permutation invariance of \(h\), reformulates the problem as a quadratic program on a portion of the sphere, and applies projected gradient descent.
The required projection, however, is tractable only in low dimensions, which limits the scalability of that approach.
These observations motivate a method that avoids sparsity guessing, scales with the ambient dimension, and recovers the entire proximal set when the operator is non-singleton.

Our approach adopts a manifold-optimization viewpoint and separates the radial and directional variables in the proximal objective.
The directional problem becomes a rank-one quadratic optimization problem on the unit sphere.
After sign and ordering reductions, a global proximal direction must arise from a local nonglobal minimizer of one of at most \(n\) sphere-constrained quadratic problems.
This observation provides both a structural description of the proximity operator and a practical route to its computation.

The main contributions are fourfold.
First, we give a complete set-valued reduction of \(\prox_{\mu h}(\vy)\), including the competition between the zero and nonzero branches.
Second, we prove that every canonical nonzero direction belongs to a finite candidate family indexed by support size and characterize each candidate by an explicit quartic equation together with positivity and second-order conditions.
Third, we develop a prefix-scan algorithm that evaluates candidate objective values without constructing all candidate vectors and that expands permutation symmetries only when the full proximal set is requested.
Fourth, we derive an exact pruning test based on the secular equation, which avoids unnecessary root solves without changing the candidate set.
Unlike the approaches in \cite{Tao2022,Jia2024}, the resulting method neither guesses the sparsity nor relies on a projection whose cost grows prohibitively with the dimension.

The remainder of the paper is organized as follows.
\Cref{sec:prox} derives
the sphere reformulation and the sign reduction.  \Cref{sec:characterization}
establishes the finite-candidate characterization and the associated quartic
equation.  \Cref{sec:algorithm} develops the direct and optimized prefix-scan
algorithms.  \Cref{sec:improvement} introduces the exact pruning test, and
\Cref{sec:numeric} reports numerical comparisons, scaling results, and
residual tests.

\section{Proximity operator via manifold optimization}\label{sec:prox}

We now reformulate the computation of \(\prox_{\mu h}(\vy)\), for \(\mu>0\) and \(\vy\in\bR^n\), as a smooth directional optimization problem on the unit sphere.
This reduction separates the radial variable, which can be minimized explicitly, from the direction, which carries the nonconvex structure of the ratio.

For every \(\vx\in\bR^n\), \(\|\vx\|_2\leq\|\vx\|_1\leq\sqrt n\|\vx\|_2\).
Since \(a\in[0,1]\), it follows immediately that \(\prox_{\mu h}(\vzero)=\{\vzero\}\).
We therefore assume throughout the remainder of this section that \(\vy\neq\vzero\).

Following the manifold-optimization viewpoint in \cite{Jia2024}, represent
every \(\vx\neq\vzero\) in polar form as
\begin{equation}\label{eq:x-polar}
  \vx=r\vu,
  \qquad r\in\bR\setminus\{0\},
  \qquad \vu\in\bS^{n-1},
\end{equation}
where \(\bS^{n-1}:=\{\vu\in\bR^n:\|\vu\|_2=1\}\).
Allowing a signed radius is convenient because the pairs \((r,\vu)\) and \((-r,-\vu)\) represent the same nonzero point.
Substituting \cref{eq:x-polar} into the
proximal objective gives
\begin{align}\label{eq:Qhy}
  Q_{h,\vy}(\vx)
  =\frac12\|\vy-\vx\|_2^2+\mu h(\vx)
  =\frac12\|\vy-r\vu\|_2^2+\mu\|\vu\|_1,
\end{align}
because \(\|r\vu\|_1/\|r\vu\|_2=\|\vu\|_1\).  Completing the square in
\(r\) yields
\begin{align}\label{eq:polar-objective}
  Q_{h,\vy}(\vx)
   & =\frac12\|\vy\|_2^2
  +\frac12\bigl(r-\langle\vy,\vu\rangle\bigr)^2
  + G_{\vy}(\vu),
\end{align}
where
\begin{align}\label{eq:G-y}
  G_{\vy}(\vu) & :=-\frac12\langle\vy,\vu\rangle^2+\mu\|\vu\|_1, \qquad \vu\in\bS^{n-1}.
\end{align}
Define
\begin{equation}\label{eq:U-y}
  \cU(\vy):=\operatorname*{argmin}_{\vu\in\bS^{n-1}}
  G_{\vy}(\vu), \qquad m(\vy):=\min_{\vu\in\bS^{n-1}}G_{\vy}(\vu).
\end{equation}
The set \(\cU(\vy)\) is nonempty because the sphere is compact and \(G_{\vy}\) is continuous.
A small but important point is that the optimal radius never vanishes at a minimizing direction.
Indeed, if \(y_j\neq0\), then \[ m(\vy) \leq G_{\vy}(\sign(y_j)\ve_j) =\mu-\tfrac12y_j^2<\mu.
\]
On the other hand, a direction satisfying
\(\langle\vy,\vu\rangle=0\) has
\(G_{\vy}(\vu)=\mu\|\vu\|_1\geq\mu\).
Hence every \(\vu\in\cU(\vy)\) has nonzero correlation with \(\vy\), and the optimal radius for that direction is \(r=\langle\vy,\vu\rangle\neq0\).
It follows that the minimizers on the
nonzero branch are precisely
\begin{equation}\label{eq:X-y}
  \cX(\vy):=
  \{\langle\vy,\vu\rangle\vu:\vu\in\cU(\vy)\}.
\end{equation}

The next proposition compares the best nonzero value with the value at the origin and thereby completes the radial reduction.
\begin{prop}\label{prop:set-valued-prox}
  For \(\vy\neq\vzero\),
  \begin{equation}\label{eq:prox-set-valued}
    \prox_{\mu h}(\vy)=
    \begin{cases}
      \{\vzero\},             & m(\vy)>\mu a, \\[1mm]
      \{\vzero\}\cup\cX(\vy), & m(\vy)=\mu a, \\[1mm]
      \cX(\vy),               & m(\vy)<\mu a.
    \end{cases}
  \end{equation}
\end{prop}
\begin{proof}
  By the preceding argument, the best value on the nonzero branch is \(\tfrac12\|\vy\|_2^2+m(\vy)\), whereas \(Q_{h,\vy}(\vzero)=\tfrac12\|\vy\|_2^2+\mu a\).
  The three cases in \cref{eq:prox-set-valued} follow by comparing these two values.
  When they are equal, both the origin and every point in \(\cX(\vy)\) are global minimizers.
\end{proof}

Thus computing \(\prox_{\mu h}(\vy)\) reduces to computing \(m(\vy)\) and \(\cX(\vy)\), or equivalently the set \(\cU(\vy)\) in \cref{eq:U-y}.
The reduced problem is still nonsmooth because the \(\ell_1\) term depends on the signs of the components of \(\vu\).
The next step removes this remaining sign dependence.

Introduce the nonnegative and positive orthants of the unit sphere:
\begin{align*}
  \bS_{+}^{n-1}
   & :=\{\vu\in\bS^{n-1}:u_i\geq0,\ 1\leq i\leq n\}, \\
  \bS_{++}^{n-1}
   & :=\{\vu\in\bS^{n-1}:u_i>0,\ 1\leq i\leq n\}.
\end{align*}
Let \(\vz:=|\vy|\), where the absolute value is taken componentwise.
On
\(\bS_{+}^{n-1}\), define
\begin{align}\label{eq:F-z}
  F_{\vz}(\vv) =-\frac12\langle\vz,\vv\rangle^2 +\mu\langle\vone_n,\vv\rangle.
\end{align}
The set of minimizers of \(F_{\vz}\) over \(\bS_{+}^{n-1}\) is denoted by \(\cV(\vz)\):
\begin{align}\label{eq:V-z}
  \cV(\vz) :=\operatorname*{argmin}_{\vv\in\bS_+^{n-1}}
  F_{\vz}(\vv).
\end{align}
For any \(\vu\in\bS^{n-1}\), \(|\langle\vy,\vu\rangle|\leq\langle|\vy|,|\vu|\rangle\), so the minimum in \cref{eq:U-y} equals the minimum in \cref{eq:V-z}.
Equality at a minimizing direction forces all nonzero products \(y_i u_i\) to have a common sign; the two possible global signs generate the same point \(\langle\vy,\vu\rangle\vu\).
Moreover, if \(z_i=0\), then every \(\vv\in\cV(\vz)\) satisfies \(v_i=0\).
Indeed, if \(v_i>0\), choose \(j\) with \(z_j>0\) and replace \((v_j,v_i)\) by \((\sqrt{v_j^2+v_i^2},0)\).
The Euclidean norm is preserved, the \(\ell_1\) term does not increase, and the correlation with \(\vz\) increases strictly, contradicting optimality.
Consequently, we have
\begin{align}
  m(\vy)
   & =\min_{\vv\in\bS_+^{n-1}}F_{\vz}(\vv),                       \label{eq:m-sign-reduction} \\
  \cX(\vy)
   & =\bigl\{\langle\vz,\vv\rangle
  (\vv\odot\sign(\vy)):\vv\in\cV(\vz)\bigr\},            \label{eq:X-sign-reduction}
\end{align}
where
\begin{equation*}
  \sign(t)=
  \begin{cases}
    1,  & t>0, \\[1mm]
    0,  & t=0, \\[1mm]
    -1, & t<0,
  \end{cases}
\end{equation*}
componentwise, and \(\odot\) denotes componentwise multiplication.
The condition \(v_i=0\) whenever \(z_i=0\) ensures that the zero entries in \(\sign(\vy)\) do not alter the norm of the recovered direction.

In summary, the original nonsmooth nonconvex problem has been reduced to the smooth quadratic problem \cref{eq:V-z} over the nonnegative orthant of the unit sphere.
Once \(\cV(\vz)\) is known, \cref{eq:X-sign-reduction} recovers the nonzero branch in the original coordinates, and \Cref{prop:set-valued-prox} determines whether the origin must also be included.
The next section gives a finite structural characterization of \(\cV(\vz)\).

\section{Structural characterization}\label{sec:characterization}
We now characterize the minimizer set \(\cV(\vz)\) in \cref{eq:V-z}.
The central observation is that the positive part of a canonical minimizer is a local nonglobal minimizer of a rank-one quadratic function on a lower-dimensional sphere.
Since there is at most one such candidate for each support size, the original continuum problem reduces to a finite comparison.

We begin by sorting the components of \(\vz\) in nonincreasing order.
Define
the monotone cone
\begin{align*}
  \bR_{\supdown}^n:=\{\vv\in\bR^n:v_1\geq\cdots\geq v_n\}
\end{align*}
and let \(\mP\) be any permutation matrix such that
\begin{align*}
  \tvy:=\mP\vz\in\bR_{\supdown}^n.
\end{align*}
Because \(\mP\) is fixed, \(\cV(\vz)=\mP^{-1}\cV(\tvy)\).
Equal components of \(\tvy\) generate a permutation symmetry.
Let
\[
  \mathfrak P(\tvy)
  :=\{\mR:\mR\text{ is a permutation matrix and }\mR\tvy=\tvy\}
\]
and define the canonical minimizer set
\begin{equation}\label{eq:canonical-set}
  \cV_{\supdown}(\tvy)
  :=\cV(\tvy)\cap\bR_{\supdown}^n.
\end{equation}

\begin{lem}\label{lem:ordering}
  If \(\vv\in\cV(\tvy)\) and \(\ty_i>\ty_j\) for \(i<j\), then \(v_i\geq v_j\).
  Moreover,
  \begin{equation}\label{eq:symmetry-orbits}
    \cV(\tvy)
    =\{\mR\vv:\vv\in\cV_{\supdown}(\tvy),\ \mR\in\mathfrak P(\tvy)\}.
  \end{equation}
\end{lem}
\begin{proof}
  If \(v_i<v_j\), interchange the two components to obtain \(\vv'\).
  The linear term is unchanged and \[ \langle\tvy,\vv'\rangle-\langle\tvy,\vv\rangle =(\ty_i-\ty_j)(v_j-v_i)>0.
  \]
  Both correlations are nonnegative, hence
  \(F_{\tvy}(\vv')<F_{\tvy}(\vv)\), a contradiction.
  Within each block of equal components of \(\tvy\), sorting the corresponding components of \(\vv\) leaves the objective unchanged.
  This yields a canonical nonincreasing representative, and all representatives in its orbit are obtained by matrices in \(\mathfrak P(\tvy)\).
\end{proof}

The lemma shows that \(\cV_{\supdown}(\tvy)\) contains all the information needed to recover \(\cV(\tvy)\).
Every vector in the canonical set is nonnegative and nonincreasing; consequently, its positive entries form a prefix.
If that prefix has length \(k\), the remaining \(n-k\) entries are zero.
Determining the prefix and its length therefore determines the entire canonical direction.

We next relate this prefix to existing results for quadratic minimization on spheres \cite{Martinez1994,Tao1995}.
For
\(1\leq k\leq n\), let \(\tvy_{:k}\) denote the first \(k\) components of
\(\tvy\) and consider
\begin{equation}\label{eq:min-F-vy-k}
  F_{\tvy_{:k}}(\vu) :=-\frac12\vu^{\mathsf T}\tvy_{:k}\tvy_{:k}^{\mathsf T}\vu +\mu\vu^{\mathsf T}\vone_k, \qquad \vu\in\bS^{k-1}.
\end{equation}
A positive local minimizer of this full-sphere problem cannot be global,
because for every \(\vu\in\bS_{++}^{k-1}\),
\begin{equation}\label{eq:positive-not-global}
  F_{\tvy_{:k}}(-\vu) =F_{\tvy_{:k}}(\vu)-2\mu\vone_k^{\mathsf T}\vu <F_{\tvy_{:k}}(\vu).
\end{equation}
This explains why local \emph{nonglobal} minimizers are the relevant positive candidates.

Set \(\cC_1:=\{1\}\), and for \(2\leq k\leq n\) define
\begin{align}\label{eq:Ck}
  \cC_k
  :=\bigl\{\vu\in\bS_{++}^{k-1}\cap\bR_{\supdown}^k:
  \vu\text{ is a local nonglobal minimizer of }
  F_{\tvy_{:k}}\text{ on }\bS^{k-1}\bigr\},
\end{align}
and
\begin{align*}\label{eq:Ck-padded}
  \bcC_k:=\{(\vu,\vzero_{n-k}):\vu\in\cC_k\}.
\end{align*}
The uniqueness result in \cite[Theorem~3.1(ii)]{Martinez1994} implies that each \(\cC_k\) is empty or a singleton.
The same holds for \(\bcC_k\).

Thus the search for \(\cV_{\supdown}(\tvy)\) reduces to the finitely many sets \(\cC_k\).
The index \(k\) has a direct interpretation: it is the support size of the canonical direction, or equivalently the length of its positive prefix.
The case \(k=1\) is always represented by \(\cC_1\), so the finite candidate union is nonempty even when none of the higher-dimensional local nonglobal minimizers exists.

\begin{thm}\label{thm:finite-candidates}
  The canonical minimizer set satisfies
  \begin{equation}\label{eq:finite-candidate-argmin}
    \cV_{\supdown}(\tvy)
    =\operatorname*{argmin}_{\vu\in\bigcup_{k=1}^n\bcC_k}
    F_{\tvy}(\vu).
  \end{equation}
  In particular, the right-hand side is a nonempty set with at most \(n\) elements.
\end{thm}
\begin{proof}
  Take \(\vv\in\cV_{\supdown}(\tvy)\).
  Since \(\vv\geq0\) is nonincreasing, it has the form \(\vv=(\vu,\vzero_{n-k})\) for some \(1\leq k\leq n\) and \(\vu\in\bS_{++}^{k-1}\).
  The point \(\vu\) lies in the relative interior of the positive orthant of \(\bS^{k-1}\).
  Global optimality of \(\vv\) over \(\bS_+^{n-1}\) therefore implies that \(\vu\) is a local minimizer of \(F_{\tvy_{:k}}\) on the full sphere.
  By \cref{eq:positive-not-global}, it is nonglobal, so \(\vv\in\bcC_k\).
  Hence \(\cV_{\supdown}(\tvy)\subseteq\bigcup_k\bcC_k\).

  Conversely, every point in \(\bigcup_k\bcC_k\) belongs to \(\bS_+^{n-1}\cap\bR_{\supdown}^n\).
  Since the union contains every canonical global minimizer, minimizing over the union gives the same value as minimizing over \(\bS_+^{n-1}\).
  The minimizers over the union are therefore exactly \(\cV_{\supdown}(\tvy)\).
\end{proof}

\Cref{thm:finite-candidates} is constructive as well as descriptive.
For each \(k\), it suffices to determine whether \(\cC_k\) is empty and, if not, to compute its unique element.
After padding the resulting prefixes with zeros, every candidate can be evaluated with the same full-dimensional objective \(F_{\tvy}\).
Retaining all candidates that attain the smallest value yields the complete canonical minimizer set, including possible ties between different support sizes.
The remaining task is therefore to characterize each \(\cC_k\) by quantities computable from the first \(k\) sorted data entries.

We next give a necessary-and-sufficient characterization of each \(\cC_k\).
For \(2\leq k\leq n\), write
\begin{equation}\label{eq:prefix-quantities}
  \begin{aligned}
    S_{1,k} & :=\|\tvy_{:k}\|_1, & \qquad S_{2,k} & :=\|\tvy_{:k}\|_2^2, \\ D_k & :=kS_{2,k}-S_{1,k}^2, & \qquad L_k & :=S_{2,k}-\ty_kS_{1,k}.
  \end{aligned}
\end{equation}
Notice that \(D_k\geq0\) by the Cauchy--Schwarz inequality and \(L_k=\sum_{i=1}^k\ty_i(\ty_i-\ty_k)\geq0\).
The quantity \(D_k\) measures the deviation of the prefix from a constant vector, while \(L_k\) will encode positivity of the last, and hence smallest, component of a candidate.
Define, for
\(0<\lambda<S_{2,k}\),
\begin{equation}\label{eq:phi-k}
  \varphi_k(\lambda)
  :=\frac{\mu^2}{S_{2,k}}
  \left[
    \frac{S_{1,k}^2}{(S_{2,k}-\lambda)^2}
    +\frac{D_k}{\lambda^2}
    \right].
\end{equation}
The associated quartic is
\begin{equation}\label{eq:psi-k}
  \psi_k(\lambda)=\sum_{i=0}^4a_i\lambda^i,
\end{equation}
where
\begin{equation}\label{eq:quartic-coefficients}
  a_0=-\mu^2S_{2,k}
  D_k, \quad a_1=2\mu^2D_k, \quad a_2=S_{2,k}^2-k\mu^2, \quad a_3=-2S_{2,k}, \quad a_4=1.
\end{equation}

\begin{thm}\label{thm:quartic-iff}
  Let \(2\leq k\leq n\).
  The following statements are equivalent.
  \begin{enumerate}[(i)]
    \item \(\cC_k\neq\emptyset\).
    \item There exists \(\lambda_k^*\in(L_k,S_{2,k})\) such that
          \begin{equation}\label{eq:root-conditions}
            \psi_k(\lambda_k^*)=0,
            \qquad
            \psi_k'(\lambda_k^*)<0.
          \end{equation}
  \end{enumerate}
  When these conditions hold, \(\lambda_k^*\) is unique and
  \(\cC_k=\{\vu_k\}\), where
  \begin{equation}\label{eq:vu-lambda-k}
    \vu_k
    =\frac{\mu}{\lambda_k^*(S_{2,k}-\lambda_k^*)}
    \left[S_{1,k}\tvy_{:k}-(S_{2,k}-\lambda_k^*)\vone_k\right].
  \end{equation}
\end{thm}
\begin{proof}
  Let \(\vz=\tvy_{:k}\), \(S_1=S_{1,k}\), \(S_2=S_{2,k}\), \(D=D_k\), and \(L=L_k\).
  For \(0<\lambda<S_2\), the stationarity
  system for the sphere-constrained problem is
  \begin{equation}\label{eq:stationarity-system}
    (\lambda\mI_k-\vz\vz^{\mathsf T})\vu+\mu\vone_k=\vzero.
  \end{equation}
  The Sherman--Morrison identity \cite{Sherman1950} gives \[ (\lambda\mI_k-\vz\vz^{\mathsf T})^{-1} =\frac1\lambda \left(\mI_k-\frac{\vz\vz^{\mathsf T}}{S_2-\lambda}\right).
  \]
  Consequently, the solution of \cref{eq:stationarity-system} is exactly the
  vector in \cref{eq:vu-lambda-k}, and its squared norm is
  \(\varphi_k(\lambda)\).  Direct expansion yields
  \begin{equation}\label{eq:phi-psi-relation}
    \varphi_k(\lambda)
    =1-\frac{\psi_k(\lambda)}{\lambda^2(S_2-\lambda)^2}.
  \end{equation}
  The smallest component of the numerator in \cref{eq:vu-lambda-k} is \(S_1\ty_k-(S_2-\lambda)=\lambda-L\).
  Thus
  \begin{equation}\label{eq:positivity-lambda}
    \vu_k\in\bS_{++}^{k-1}
    \quad\Longleftrightarrow\quad
    \lambda>L.
  \end{equation}
  The formula also shows that \(\vu_k\) is nonincreasing.

  Suppose first that \(\vu\in\cC_k\).
  The local-nonglobal characterization in \cite[Theorem~3.1(ii)]{Martinez1994} associates with \(\vu\) the larger simple root of the secular equation \(\varphi_k(\lambda)=1\) in \((0,S_2)\).
  Hence \(\varphi_k'(\lambda)>0\).
  By \cref{eq:phi-psi-relation}, at a root \[ \varphi_k'(\lambda) =-\frac{\psi_k'(\lambda)}{\lambda^2(S_2-\lambda)^2}, \] so \(\psi_k'(\lambda)<0\).
  Positivity gives \(\lambda>L\), proving necessity and \cref{eq:vu-lambda-k}.

  Conversely, suppose \cref{eq:root-conditions} holds and define \(\vu_k\) by \cref{eq:vu-lambda-k}.
  Equations \cref{eq:phi-psi-relation,eq:positivity-lambda} show that \(\|\vu_k\|_2=1\) and \(\vu_k\in\bS_{++}^{k-1}\cap\bR_{\supdown}^k\).
  It also satisfies \cref{eq:stationarity-system}.
  Let \(\alpha:=\vz^{\mathsf T}\vu_k\).
  For every tangent vector
  \(\vv\perp\vu_k\),
  \begin{align}
    \vv^{\mathsf T}(\lambda_k^*\mI_k-\vz\vz^{\mathsf T})\vv
     & \geq
    \bigl(\lambda_k^*-S_2+\alpha^2\bigr)\|\vv\|_2^2.              \label{eq:tangent-bound}
  \end{align}
  To see this bound, note that the orthogonal projection of \(\vz\) onto the tangent space at \(\vu_k\) has squared norm \(S_2-\alpha^2\).
  Cauchy--Schwarz therefore gives \((\vz^{\mathsf T}\vv)^2\leq(S_2-\alpha^2)\|\vv\|_2^2\), which is equivalent to \cref{eq:tangent-bound}.
  At a secular root, a direct calculation gives
  \begin{equation}\label{eq:phi-prime-second-order}
    \varphi_k'(\lambda_k^*)
    =\frac{2}{\lambda_k^*(S_2-\lambda_k^*)}
    \bigl(\lambda_k^*-S_2+\alpha^2\bigr).
  \end{equation}
  Because \(\psi_k'(\lambda_k^*)<0\), the left-hand side is positive.
  Therefore the Lagrangian Hessian is positive definite on the tangent space, and \(\vu_k\) is a strict local minimizer.
  It is not global by \cref{eq:positive-not-global}; hence \(\vu_k\in\cC_k\).

  Finally, if \(D>0\), then \[ \varphi_k''(\lambda) =\frac{6\mu^2}{S_2} \left[\frac{S_1^2}{(S_2-\lambda)^4}+\frac{D}{\lambda^4}\right]>0, \] so there is at most one secular root with positive derivative.
  If \(D=0\), \(\varphi_k\) is strictly increasing.
  This proves uniqueness.
\end{proof}

The quartic \(\psi_k\) has real coefficients, and its roots admit closed-form expressions \cite{Chavez-Pichardo2022}.
Those expressions are lengthy and do not clarify the structural conditions in \Cref{thm:quartic-iff}; in numerical work, a stable polynomial or bracketed secular-equation solver is preferable.
The following example nevertheless gives an explicit calculation in the important degenerate case in which the entries of \(\tvy_{:k}\) are equal.

\begin{example}\label{ex:equal-components}
  Let \(k\geq2\) and \(\tvy_{:k}=y\vone_k\) with \(y>0\).
  Then \[ \psi_k(\lambda) =\lambda^2\bigl((\lambda-ky^2)^2-k\mu^2\bigr).
  \]
  The only admissible root is
  \(\lambda_k^*=ky^2-\sqrt{k}\mu\), and it satisfies
  \cref{eq:root-conditions} exactly when
  \[
    \mu<y^2\sqrt{k}.
  \]
  In that case \(\cC_k=\{k^{-1/2}\vone_k\}\); otherwise
  \(\cC_k=\emptyset\).
  This statement concerns only \(k\geq2\).
  The separately treated set \(\cC_1=\{1\}\) is always present, so the finite candidate union is never empty.
\end{example}

\section{Algorithms and complexity}\label{sec:algorithm}

\Cref{thm:finite-candidates,thm:quartic-iff} yield an exact finite procedure
for computing the directional minimizers.
For each support size \(k\), \Cref{thm:quartic-iff} determines whether \(\cC_k\) is empty and constructs its unique element when it is not.
\Cref{thm:finite-candidates} then selects
every padded candidate that attains the smallest value of \(F_{\tvy}\).
Keeping all ties is essential because the proximity operator can be set-valued \cite{Tao2022}.
After the canonical minimizers are identified,
\cref{eq:symmetry-orbits,eq:X-sign-reduction} recover the complete nonzero
candidate set in the original coordinates:
\begin{equation}\label{eq:X-final}
  \cX(\vy)
  =\left\{
  \langle\tvy,\vu\rangle
  \bigl(\mP^{-1}\mR\vu\odot\sign(\vy)\bigr):
  \vu\in\cV_{\supdown}(\tvy),\ \mR\in\mathfrak P(\tvy)
  \right\}.
\end{equation}
Accordingly,
\begin{equation}\label{eq:prox-h-vy-final}
  \prox_{\mu h}(\vy)=
  \begin{cases}
    \{\vzero\},             & m(\vy)>\mu a, \\[1mm]
    \{\vzero\}\cup\cX(\vy), & m(\vy)=\mu a, \\[1mm]
    \cX(\vy),               & m(\vy)<\mu a.
  \end{cases}
\end{equation}
Thus the finite search returns all canonical proximal directions, while the explicit symmetry expansion returns the full proximal set.
If only one proximal point is required, one representative from each tie orbit need not be constructed; enumerating the entire set is naturally output-sensitive.

We summarize the above procedure in \Cref{alg:prox-h-naive}.

\begin{algorithm}[H]
  \caption{Finite-candidate computation of \(\prox_{\mu h}(\vy)\) (direct version).}
  \label{alg:prox-h-naive}
  \begin{algorithmic}[1]
    \Statex \textbf{Input:} \(\vy\in\bR^n\), \(\mu>0\), and \(a\in[0,1]\)
    \Statex \textbf{Output:} \(\prox_{\mu h}(\vy)\)
    \If{\(\vy=\vzero\)}
    \State \Return \(\{\vzero\}\)
    \EndIf
    \State Find \(\mP\) such that \(\tvy=\mP|\vy|\in\bR_{\supdown}^n\).
    \State Set \(\cC_1=\{1\}\) and \(\bcC_1=\{(1,\vzero_{n-1})\}\).
    \For{\(k=2,\ldots,n\)}
    \State Form \(S_{1,k},S_{2,k},D_k,L_k\) and the coefficients in \cref{eq:quartic-coefficients}.
    \State Compute the real roots of \(\psi_k\).
    \If{there is a root \(\lambda_k^*\in(L_k,S_{2,k})\) with \(\psi_k'(\lambda_k^*)<0\)}
    \State Form \(\vu_k\) by \cref{eq:vu-lambda-k} and set
    \(\bcC_k=\{(\vu_k,\vzero_{n-k})\}\).
    \Else
    \State Set \(\bcC_k=\emptyset\).
    \EndIf
    \EndFor
    \State Compute \(\cV_{\supdown}(\tvy)\) by comparing \(F_{\tvy}\) on
    \(\bigcup_{k=1}^n\bcC_k\).
    \State Expand equal-component symmetries as in \cref{eq:X-final} and apply \cref{eq:prox-h-vy-final}.
  \end{algorithmic}
\end{algorithm}

If prefix quantities, candidate vectors, and objective values are recomputed from scratch, \Cref{alg:prox-h-naive} costs \(O(n^2)\) after sorting.
This cost is avoidable because the quantities in \cref{eq:prefix-quantities} admit constant-time updates, and candidate-vector construction can be deferred until the winning support sizes are known.
For a candidate root, substituting \cref{eq:vu-lambda-k} gives
\begin{equation}\label{eq:F-candidate-scalar}
  f_k
  :=F_{\tvy}(\vu_k,\vzero_{n-k})
  =-\frac{\mu^2S_{1,k}^2}{2(S_{2,k}-\lambda_k^*)^2}
  +\frac{\mu^2\bigl(S_{1,k}^2-k(S_{2,k}-\lambda_k^*)\bigr)}
  {\lambda_k^*(S_{2,k}-\lambda_k^*)}.
\end{equation}
This is an \(O(1)\) calculation once the prefix quantities and root are available.
We can therefore maintain the two prefix sums and the current minimum objective value while scanning through the support sizes, without forming a length-\(k\) vector at every iteration.
This observation yields the optimized procedure in \Cref{alg:prox-h-optimized}.

\begin{algorithm}[H]
  \caption{Scalar prefix scan for \(\prox_{\mu h}(\vy)\) (optimized version).}
  \label{alg:prox-h-optimized}
  \begin{algorithmic}[1]
    \Statex \textbf{Input:} \(\vy\in\bR^n\), \(\mu>0\), and \(a\in[0,1]\)
    \Statex \textbf{Output:} \(\prox_{\mu h}(\vy)\), where
    \(h(\vx)=\|\vx\|_1/\|\vx\|_2\) for \(\vx\neq\vzero\)
    \If{\(\vy=\vzero\)}
    \State \Return \(\{\vzero\}\)
    \EndIf
    \State Sort: \(\tvy=\mP|\vy|\in\bR_{\supdown}^n\).
    \State Initialize
    \(S_1=\ty_1\), \(S_2=\ty_1^2\),
    \(f_{\min}=-\tfrac12\ty_1^2+\mu\), and
    \(\mathcal I=\{(1,\varnothing)\}\).
    \For{\(k=2,\ldots,n\)}
    \State Update \(S_1\gets S_1+\ty_k\) and \(S_2\gets S_2+\ty_k^2\).
    \State Set \(D=kS_2-S_1^2\), \(L=S_2-\ty_kS_1\), and form \(\psi_k\).
    \State Find the unique root \(\lambda_k^*\in(L,S_2)\) with
    \(\psi_k(\lambda_k^*)=0\) and \(\psi_k'(\lambda_k^*)<0\), if it exists.
    \If{such a root exists}
    \State Compute \(f_k\) from \cref{eq:F-candidate-scalar}.
    \If{\(f_k<f_{\min}\)}
    \State \(f_{\min}\gets f_k\), \(\mathcal I\gets\{(k,\lambda_k^*)\}\).
    \ElsIf{\(f_k=f_{\min}\)}
    \State \(\mathcal I\gets\mathcal I\cup\{(k,\lambda_k^*)\}\).
    \EndIf
    \EndIf
    \EndFor
    \State Set \(m(\vy)\gets f_{\min}\).
    \State For each record in \(\mathcal I\), construct the winning vector only once using \cref{eq:vu-lambda-k}.
    \State Expand the symmetry group and apply \cref{eq:prox-h-vy-final}.
  \end{algorithmic}
\end{algorithm}

The sorting step costs \(O(n\log n)\).
In the fixed-precision arithmetic model, prefix updates, quartic coefficients, root tests for a fixed-degree polynomial, and scalar objective comparisons cost \(O(1)\) per support size, so the post-sorting scan costs \(O(n)\).
Constructing one selected direction costs \(O(n)\).
Therefore the overall cost for returning one proximal point is \(O(n\log n)\), or \(O(n)\) when the magnitudes are already sorted.
If \(p\) distinct canonical minimizers are tied, their explicit construction costs \(O(pn)\).
Enumerating every permutation generated by equal data components is necessarily output-sensitive and can be combinatorial in the number of reported points.
The theoretical procedure is exact because it compares the complete finite candidate family.
A floating-point implementation approximates the quartic or secular roots, so the residual and tie tolerances used in practice are reported explicitly in \Cref{sec:numeric}.

\section{Further improvements with a pruning test}\label{sec:improvement}
The finite characterization reduces the directional problem from an infinite feasible set to at most \(n\) candidates.
Nevertheless, a direct implementation may still attempt a quartic root computation for every support size.
Many of these root solves are unnecessary because the corresponding sphere problem either has no local nonglobal minimizer or has one whose last component is nonpositive.
We now combine these two requirements into a single scalar test.

Specifically, we derive a necessary-and-sufficient condition for the existence of an element of \(\cC_k\).
The condition first identifies the minimum of the strictly convex secular function and then incorporates the lower bound \(\lambda>L_k\) required by the positive orthant.
Consequently, the pruned variant preserves the exact candidate family while skipping every root solve that cannot yield an admissible direction.

\begin{thm}\label{thm:prune}
  Fix \(2\leq k\leq n\) and use the quantities in \cref{eq:prefix-quantities}.
  \begin{enumerate}[(i)]
    \item If \(\ty_k=0\), then \(\cC_k=\emptyset\).
    \item If \(D_k=0\), then all components of \(\tvy_{:k}\) are equal and
          \begin{equation}\label{eq:equal-existence}
            \cC_k\neq\emptyset
            \quad\Longleftrightarrow\quad
            \mu<\ty_1^2\sqrt{k}.
          \end{equation}
    \item Suppose \(D_k>0\) and \(\ty_k>0\).  Define
          \begin{equation}\label{eq:lambda0-M}
            \lambda_{0,k}
            :=\frac{S_{2,k}}
            {1+\left(S_{1,k}^2/D_k\right)^{1/3}},
            \qquad
            M_k:=\max\{L_k,\lambda_{0,k}\}.
          \end{equation}
          Then
          \begin{equation}\label{eq:exact-pruning-test}
            \cC_k\neq\emptyset
            \quad\Longleftrightarrow\quad
            \varphi_k(M_k)<1.
          \end{equation}
  \end{enumerate}
\end{thm}
\begin{proof}
  If \(\ty_k=0\), then \(L_k=S_{2,k}\), so the positivity interval \((L_k,S_{2,k})\) is empty.
  If \(D_k=0\), equality holds in the Cauchy--Schwarz inequality, and the assertion follows from \Cref{ex:equal-components}.

  Assume \(D_k>0\) and \(\ty_k>0\).
  From \cref{eq:phi-k}, \(\varphi_k\) tends to \(+\infty\) at both endpoints of \((0,S_{2,k})\), and \[ \varphi_k''(\lambda) =\frac{6\mu^2}{S_{2,k}} \left[ \frac{S_{1,k}^2}{(S_{2,k}-\lambda)^4} +\frac{D_k}{\lambda^4} \right]>0.
  \]
  Thus \(\varphi_k\) is strictly convex and its unique minimizer is
  \(\lambda_{0,k}\).
  A local nonglobal sphere candidate corresponds to the right root \(\lambda_+>\lambda_{0,k}\) of \(\varphi_k(\lambda)=1\).
  By \Cref{thm:quartic-iff}, this candidate is strictly positive exactly when \(\lambda_+>L_k\).

  If \(L_k\leq\lambda_{0,k}\), the positivity inequality is automatic once the two simple roots exist, and existence is equivalent to \(\varphi_k(\lambda_{0,k})<1\).
  If \(L_k>\lambda_{0,k}\), the function is strictly increasing on \([L_k,S_{2,k})\), so \(\lambda_+>L_k\) is equivalent to \(\varphi_k(L_k)<1\).
  Both cases are exactly \cref{eq:exact-pruning-test}.
  The inequality is strict.
  If \(M_k=\lambda_{0,k}\), equality produces a double secular root with zero derivative.
  If \(M_k=L_k>\lambda_{0,k}\), equality places the increasing root at \(\lambda=L_k\), where the last component in \cref{eq:vu-lambda-k} is zero.
  Neither case yields an element of \(\cC_k\).
\end{proof}

\begin{example}\label{eg:positive-threshold}
  The following two-dimensional instance illustrates why the positivity threshold \(L_k\) must be incorporated into the pruning criterion.
  Let \[ \tvy=(2,1)^{\mathsf T} \qquad\text{and}\qquad \mu=2.03.
  \]
  For \(k=2\), we have
  \[
    S_{1,2}=3,\qquad S_{2,2}=5,\qquad
    D_2=1,\qquad L_2=2,
  \]
  and
  \[
    \lambda_{0,2}
    =\frac{5}{1+9^{1/3}}
    \approx 1.6233.
  \]
  Although
  \[
    \varphi_2(\lambda_{0,2})\approx 0.9633<1,
  \]
  we have \(L_2>\lambda_{0,2}\), and hence
  \(M_2=\max{L_2,\lambda_{0,2}}=L_2\).
  Moreover, \[ \varphi_2(M_2)=\varphi_2(L_2)\approx 1.0302>1.
  \]
  Therefore, Theorem~\ref{thm:prune} gives
  \(\mathcal C_2=\varnothing\).
  This example shows that the full-sphere condition \(\varphi_k(\lambda_{0,k})<1\) alone does not guarantee a positive candidate, whereas the corrected criterion \(\varphi_k(M_k)<1\) correctly incorporates the positive-orthant requirement.
\end{example}

The test in \Cref{thm:prune} uses only \(S_{1,k}\), \(S_{2,k}\), \(\ty_k\), and \(k\), and therefore costs \(O(1)\) per support size.
It can be inserted immediately before the root computation in \Cref{alg:prox-h-optimized}.
When the test succeeds, the desired root is the unique solution of \(\varphi_k(\lambda)=1\) on the increasing interval \((M_k,S_{2,k})\).
A bracketed monotone solver can therefore be used without computing all four quartic roots.
When the test fails, \(\cC_k\) is empty.
The pruning test does not change the \(O(n)\) post-sorting complexity, but it can substantially reduce the number of comparatively expensive root solves, as quantified in \Cref{sec:numeric}.

\section{Numerical experiments}\label{sec:numeric}
We present numerical experiments that compare the proposed method with the algorithm of \cite{Tao2022} and assess scaling, pruning effectiveness, and floating-point residuals.
To match the baseline, we set \(a=1\) in \cref{def:h}.
The finite characterization identifies the global proximal solution, whereas the sparsity-selection heuristic of \cite{Tao2022} can return a suboptimal point on particular inputs.

Because the method in \cite{Tao2022} works directly with the original proximal objective \(Q_{h,\vy}\) in \cref{eq:Qhy}, we report that objective for both methods.
By definition, \(\prox_{\mu h}(\vy)\) is the set of global minimizers of \(Q_{h,\vy}\).
The baseline values were obtained from the implementation accompanying \cite{Tao2022}.
Our comparisons concern the returned points for the stated inputs and parameter values; the exactness of the proposed method follows from the finite-candidate theory rather than from a numerical comparison.

For the quartic implementation, we scale the variable as \(\lambda=S_{2,k}x\), compute the polynomial roots in \(x\), retain a root as real when its imaginary part is at most \(10^{-8}\max\{1,|\operatorname{Re}x|\}\), and then apply the interval and derivative tests in \Cref{thm:quartic-iff}.
In the pruned implementation, after \Cref{thm:prune} succeeds, we instead solve \(\varphi_k(\lambda)=1\) on the increasing interval \((M_k,S_{2,k})\) by Brent's method, with scale-aware absolute tolerance \(10^{-13}\max\{1,S_{2,k}\}\) and relative tolerance \(2\times10^{-14}\).
Objective ties are tested with a scaled tolerance of \(10^{-10}\).
For every accepted root we monitor the secular residual \(|\varphi_k(\lambda_k^*)-1|\), the norm residual \(|\|\vu_k\|_2-1|\), and the stationarity residual \[ \left\| (\lambda_k^*\mI_k-\tvy_{:k}\tvy_{:k}^{\mathsf T})\vu_k +\mu\vone_k \right\|_2.
\]

\subsection{Comparison on the two illustrative examples}

\begin{example}\label{eg:comparison-2}
  Let \(\vy=(4,4,3,3,2,2)^{\mathsf T}\).
  The complete candidate comparison produces the results in \Cref{tab:comparison-2}.
  \begin{table}[H]
    \centering
    \small
    \resizebox{\textwidth}{!}{%
      \begin{tabular}{c c c c c}
        \toprule
        \(\mu\) & Method         & Selected \(k\) & Returned point                                        & \(Q_{h,\vy}\) \\
        \midrule
        1       & Proposed       & 6              & \([4.033,4.033,2.990,2.990,1.948,1.948]^{\mathsf T}\) & 2.360         \\
        1       & \cite{Tao2022} & 6              & \([4.033,4.033,2.990,2.990,1.948,1.948]^{\mathsf T}\) & 2.360         \\
        13      & Proposed       & 6              & \([4.455,4.455,2.423,2.423,0.392,0.392]^{\mathsf T}\) & 29.403        \\
        13      & \cite{Tao2022} & 3              & \([3.597,3.597,1.542,0,0,0]^{\mathsf T}\)             & 31.091        \\
        \bottomrule
      \end{tabular}%
    }
    \caption{Objective comparison for \Cref{eg:comparison-2}.}
    \label{tab:comparison-2}
  \end{table}
  For \(\mu=1\), the two procedures agree.
  For \(\mu=13\), the complete candidate comparison selects \(k=6\) and attains a smaller objective than the output of \cite{Tao2022}.
  Even at \(k=3\), the candidate generated by \cref{eq:vu-lambda-k} has objective \(30.610\), below \(31.091\) for the reported baseline point.
  These computations illustrate that the sparsity-selection heuristic in \cite{Tao2022} can return a suboptimal point; they are not used as a proof of exactness.
\end{example}

\begin{example}\label{eg:comparison-1}
  Let \(\vy=(9,7,6,4,2)^{\mathsf T}\).
  The results are shown in \Cref{tab:comparison-1}.
  \begin{table}[H]
    \centering
    \small
    \resizebox{\textwidth}{!}{%
      \begin{tabular}{c c c c c}
        \toprule
        \(\mu\) & Method         & Selected \(k\) & Returned point                                  & \(Q_{h,\vy}\) \\
        \midrule
        1       & Proposed       & 5              & \([9.026,7.004,5.993,3.970,1.948]^{\mathsf T}\) & 2.051         \\
        1       & \cite{Tao2022} & 4              & \([9.021,7.000,5.989,3.968,0]^{\mathsf T}\)     & 3.926         \\
        48      & Proposed       & 4              & \([10.255,6.362,4.415,0.521,0]^{\mathsf T}\)    & 90.740        \\
        48      & \cite{Tao2022} & 2              & \([8.506,4.642,0,0,0]^{\mathsf T}\)             & 96.030        \\
        \bottomrule
      \end{tabular}%
    }
    \caption{Objective comparison for \Cref{eg:comparison-1}.}
    \label{tab:comparison-1}
  \end{table}
  For both values of \(\mu\), the complete candidate comparison improves on the reported baseline objective.
  At \(\mu=48\), the candidate from our characterization at the same support size \(k=2\) has objective \(94.789\), again below the baseline value \(96.030\).
  Together with \Cref{eg:comparison-2}, this example separates two possible sources of suboptimality: selecting an incorrect support size and returning a nonoptimal point even after a support size has been selected.
  The finite candidate comparison addresses both issues by constructing the unique admissible sphere candidate for every support size and then comparing their objective values globally.
\end{example}

\subsection{Scaling, pruning, and residuals}

We next assess the post-sorting implementations on larger random problems.
For each dimension, we generated ten independent vectors with half-normal entries using random seed 20260825 and set \(\mu=0.6\sqrt n\).
The reported times are medians over the ten instances and include sorting.
The direct method recomputes prefix quantities and constructs each candidate vector; it was omitted beyond \(n=2000\).
The prefix-quartic scan uses scaled polynomial roots, while the pruned secular scan first applies \Cref{thm:prune} and then uses the bracketed secular solve described above.
The computations used Python~3.13.5, NumPy~2.3.5, and SciPy~1.17.0 on an Intel Xeon Platinum 8370C processor at 2.80 GHz.

\begin{table}[H]
  \centering
  \small
  \resizebox{\textwidth}{!}{%
    \begin{tabular}{r r r r r r r}
      \toprule
      \(n\)  & Direct (ms)           & Prefix quartic (ms) & Pruned secular (ms)
             & Full-sphere pass (\%) & Positive pass (\%)  & Max. objective gap                                        \\
      \midrule
      100    & 5.776                 & 4.426               & 1.219               & 83.3 & 54.5 & \(3.3\times10^{-12}\) \\
      250    & 14.307                & 10.567              & 2.316               & 67.2 & 46.7 & \(3.1\times10^{-12}\) \\
      500    & 30.096                & 21.765              & 5.049               & 76.9 & 51.7 & \(8.1\times10^{-12}\) \\
      1,000  & 58.033                & 42.188              & 9.519               & 75.2 & 50.7 & \(1.5\times10^{-11}\) \\
      2,000  & 121.728               & 88.651              & 19.340              & 74.0 & 50.1 & \(3.9\times10^{-11}\) \\
      5,000  & --                    & 224.149             & 51.286              & 75.0 & 50.5 & \(5.0\times10^{-11}\) \\
      10,000 & --                    & 437.637             & 103.797             & 75.4 & 50.6 & \(1.7\times10^{-10}\) \\
      \bottomrule
    \end{tabular}%
  }
  \caption{Runtime and pruning results.  ``Full-sphere pass'' is the fraction
    of support sizes satisfying \(\varphi_k(\lambda_{0,k})<1\); ``positive
    pass'' is the fraction also satisfying the exact positive-orthant test.
    The objective gap compares the prefix-quartic and pruned-secular scans.
  }
  \label{tab:runtime-random}
\end{table}

The timing data are consistent with linear post-sorting growth.
On these instances, the positive-orthant condition rejects roughly one third of the support sizes that pass the full-sphere test, thereby avoiding many unnecessary root solves.
The total asymptotic cost remains \(O(n\log n)\) because of the initial sort.
The pruned secular implementation is also consistently faster than the prefix-quartic implementation in this experiment.
This improvement reflects both the lower pass rate and the use of a bracketed solve on an interval where the secular function is monotone; it is a practical speedup rather than a change in the asymptotic order.

We also tested 100 independent problems from each of four input families at \(n=200\), using seed 20260824 and \(\mu=0.05\|\vy\|_\infty^2\sqrt n\).
The sparse vectors had approximately 70\% zero entries, and the repeated-entry vectors were obtained by rounding half-normal samples to one decimal place.
The stationarity residual was normalized by \(1+\mu\sqrt{k}\); explicitly, the reported relative residual is \[ \frac{\|\,(\lambda_k^*\mI_k- \tvy_{:k}\tvy_{:k}^{\mathsf T})\vu_k+\mu\vone_k\,\|_2} {1+\mu\sqrt{k}}.
\]
The last column is the absolute gap between the best directional objective
values produced by the scaled-quartic and bracketed-secular scans.

\begin{table}[H]
  \centering
  \small
  \resizebox{\textwidth}{!}{%
    \begin{tabular}{l r r r r}
      \toprule
      Family           & Max. secular residual               & Max. norm residual
                       & Max. relative stationarity residual & Max. objective gap    \\
      \midrule
      Half-normal      & \(3.0\times10^{-13}\)               & \(1.5\times10^{-13}\)
                       & \(8.1\times10^{-15}\)               & \(1.1\times10^{-11}\) \\
      Uniform          & \(5.3\times10^{-13}\)               & \(2.7\times10^{-13}\)
                       & \(1.2\times10^{-14}\)               & \(1.2\times10^{-11}\) \\
      Sparse           & \(1.9\times10^{-13}\)               & \(9.1\times10^{-14}\)
                       & \(2.5\times10^{-15}\)               & \(2.0\times10^{-12}\) \\
      Repeated entries & \(2.8\times10^{-13}\)               & \(1.4\times10^{-13}\)
                       & \(1.2\times10^{-14}\)               & \(1.4\times10^{-11}\) \\
      \bottomrule
    \end{tabular}%
  }
  \caption{Worst residuals and agreement over 100 random problems per family.
    These tests assess floating-point behavior; exactness follows from the theoretical characterization.
  }
  \label{tab:residual-random}
\end{table}

The uniformly small residuals show that both proposed implementations solve their secular equations accurately across all four input families.
The small objective gaps also confirm close agreement between the scaled-quartic and bracketed-secular scans.
These experiments assess numerical stability; the global characterization itself follows from the preceding theory.
The sparse family tests prefixes with trailing zeros, whereas the repeated- entry family activates the permutation symmetries described in \Cref{lem:ordering}.
Their residuals are comparable with those for continuous random inputs, indicating that these structural cases do not introduce a noticeable loss of numerical accuracy.

\section{Conclusion}
We developed a theoretically exact and numerically accurate procedure for computing the complete proximity operator of the \(\ell_1/\ell_2\) ratio.
The radial reduction separates the possible minimizer at the origin from a directional problem on the sphere.
After sign and ordering reductions, that directional problem becomes a finite collection of rank-one quadratic problems indexed by support size.
Each nontrivial candidate is generated by a unique simple quartic root satisfying both the positive-orthant inequality and the second-order derivative condition.
Comparing all admissible candidates recovers every canonical minimizer, and permutation expansion recovers the complete proximal set when ties occur.

With prefix sums, candidate testing and objective comparison form an \(O(n)\) scan after \(O(n\log n)\) sorting.
The exact pruning criterion eliminates support sizes that cannot yield a positive candidate and provides a monotone interval for the remaining secular root solve.
Explicit enumeration of tied outputs is necessarily output-sensitive.
The numerical experiments demonstrate objective improvements over a sparsity-guessing baseline, linear post-sorting growth, effective pruning, and small secular, norm, and stationarity residuals.
Here exactness refers to the finite theoretical characterization; the reported residuals quantify the accuracy of its floating-point implementation.

\section*{Acknowledgements}
We are grateful to Dr.~Min Tao for generously sharing the code accompanying \cite{Tao2022}, which facilitated our numerical comparisons.
We also thank the anonymous reviewers for their constructive comments and suggestions, which helped us improve the manuscript.
The \Cref{eg:positive-threshold} was suggested by one of the reviewers.
The work of L.~Shen was supported in part by the National Science Foundation under grant DMS-2208385, and the work of G.~Song was supported in part by the National Science Foundation under grant DMS-2318781.

\bibliographystyle{plain}
\bibliography{l1-l2-minimization_revised}

\end{document}